\documentclass[12pt,oneside,leqno]{amsart}
\usepackage[all]{xy}
\usepackage{amsfonts,amsmath,oldgerm,amssymb,amscd,mathrsfs,mathtools,amsthm}
\usepackage{enumerate}
\usepackage{xcolor}
\UseComputerModernTips
\numberwithin{equation}{section}
\usepackage[breaklinks]{hyperref}

\usepackage{float}

\makeatletter
\@namedef{subjclassname@2010}{%
	\textup{2010} Mathematics Subject Classification}
\makeatother

\usepackage{caption} 
\newtheorem{theorem}{Theorem}%[section]
\newtheorem{lemma}{Lemma}%[section]
\newtheorem{prop}{Proposition}%[section]
\newtheorem{coro}{Corollary}%[section]

\theoremstyle{definition}
\newtheorem*{question}{Question}
\newtheorem{rem}{Remark}
\newtheorem*{defn}{Definition} %[section]
\newtheorem*{conj}{Conjecture}%[section]
\newtheorem{exmp}{Example}[section]

\newenvironment{pff}{\hspace*{-\parindent}{\bf Proof \,}}
{\hfill $\Box$ \vspace*{0.2cm}}

\title{On dynamical irreducible set of polynomials}
\author{Himanshu Sharma} 
\address{\endgraf Department of Mathematics, \endgraf Tel Aviv University, Israel}
\email{himanshusharma985@gmail.com}

\usepackage{fancyhdr}
\begin{document}
	
	\date{\today}

	\begin{abstract}
	 In this article, a necessary and sufficient condition is proved for the dynamical irreducibility of a family of polynomials over a finite field. Using this result, an explicit construction of a dynamically irreducible set of polynomials is given over the finite field $\mathbb{F}_{p^p}$ where $p$ is a prime. Moreover, the existence of dynamically irreducible sets of size at least $p^2$ is also established over every finite field $\mathbb{F}_q$ where $q$ is a $p$ power. Finally, a bound on the cardinality of the set is given that needs to be tested for the dynamical irreducibility of a set of polynomials.
		 \end{abstract}
         \maketitle
	\section{Introduction}
	Let $K$ be a field, and let $f(x)\in{{K[x]}}$. Then, for $n\in \mathbb{N}\cup \{0\},$ the $n$-th iterate of $f(x)$ is defined inductively as follows: $$f^0(x)=x, \quad f^n(x)=f(f^{n-1}(x)).$$
		The polynomial $f(x)$ is said to be \textit{stable} over $K$ if $f^n(x)$ is irreducible over $K$ for each $n\geq{1}$. As an example, all the iterates of $f(x)=x^2+5$ are irreducible over $\mathbb{Q}$, that is, $x^2+5$ is stable over $\mathbb{Q}.$ Similarly, one can think of the situation where we compose more than one polynomial. We have the following definition:
        \begin{defn}
            Let $f_1,f_2,\cdots,f_n$ be irreducible polynomials over a field $K$. Then, they are said to be dynamically irreducible if all polynomials formed by the compositions of $f_1,f_2,\cdots,f_n$ are irreducible over $K$.
        \end{defn}
        Often, authors use the term \textit      {stable} when $n=1$ in the above definition. But dynamically irreducible fits in the general context. We are interested in the dynamical irreducibility of a set of polynomials  over finite fields. For $n=1$, Jones et. al. \cite{Jones} gave a criterion for the stability of quadratic polynomials over a finite field, and recently Fernandes et al. gave a nice criterion for the stability of a binomial in \cite{Panario2025}, which gives a complete understanding of the stability of binomials. \\
In the literature, we have many articles that deal with the stability of a polynomial (dynamical irreducibility of a set of polynomials with cardinality $1$) over finite fields and over number fields \cite{Jones,Fein,Sharma,Panario2025}. The focus of this article is to construct large dynamically irreducible sets of polynomials over finite fields. In particular, this article investigates the following question:
\begin{question}
    Can we construct a dynamically irreducible set of polynomials with cardinality more than 1 over any finite field.
\end{question}
We recall the following definition of $m$-free elements in a finite field.
\begin{defn}[\cite{Mullen}, Definition 5.1]
    Let $p$ be a prime and let $q$ be a $p$-power. We say that \(\alpha \in {\mathbb{F}}_q^*\) is \( m \)-free for \( m \mid q-1 \) if the equality \(\alpha = \beta^d\) with \(\beta \in \mathbb{F}_q\), for any divisor \( d \) of \( m \), implies \( d = 1 \).
\end{defn}

First, we prove the following result, which is a generalization of the criterion given by Fernandes et al. \cite{Panario2025} and gives a necessary and sufficient condition for the dynamical irreducibility of a set of certain polynomials with cardinality more than $1$.

\begin{theorem}\label{Main1}
    Let $p$ be a prime and let $q=p^r$ for some $r\geq{1}$. Assume that $f_i(x)=(x-a_i)^t+b_i$ is an irreducible polynomial over $\mathbb{F}_q$ for each $1\leq{i}\leq{k}$. The set $A=\{f_i(x):1\leq{i}\leq{k}\}$ is dynamically irreducible over $\mathbb{F}_q$ if and only if \begin{equation}\label{dynamical irr}
        f_{j_1}\circ f_{j_2}\circ f_{j_3}\circ\cdots\circ f_{j_r}(b_i)
    \end{equation} is rad(t)-free for all $1\leq{j_1,j_2,\ldots,j_r,i}\leq{k}.$
     
\end{theorem}

Next, using the above criterion, we construct a dynamically irreducible set of polynomials over a finite field $\mathbb{F}_{p^p}$ for every odd prime $p$.
\begin{theorem}\label{Main2}
    Let $p$ be an odd prime, and $ h\in{\mathbb{F}_p^{*}}$. Assume that $\beta$ is a root of the irreducible polynomial $f(x)=x^p-x-h\in{\mathbb{F}_p[x]}$. Let $f_{a,b}(x)=(x-a-\beta)^t+b+\beta\in{\mathbb{F}_{p^p}[x]}$ where $a,b\in{\mathbb{F}_p}$ and $t\geq{3}$. Then, the set 
    $$A=\{f_{a,b}(x):a,b\in{\mathbb{F}_p}\}$$
    is dynamically irreducible if and only if \begin{itemize}
        \item [(i)]$h$ and $-h$ are $rad(t)$-free,
        \item [(ii)] $rad(t)\mid p^p-1$,
        \item [(iii)] if $4\mid t$ then $4\mid p^p-1$.
    \end{itemize}
\end{theorem}
By using estimates of character sums, we generalize the above result to an arbitrary finite field and prove the following result:
\begin{theorem}\label{lower-bound}
    Let $q=p^r$ and let $t=s^k$ be an odd prime power such that $s\mid q-1$. Then, there exists an $\beta\in{\mathbb{F}_q}\setminus \mathbb{F}_p$ such that the set $A=\{(x-a-\beta)^t+b+\beta:a,b\in{\mathbb{F}_p}\}$ is dynamically irreducible. Moreover, for any odd $t\geq{3}$ with $rad(t)\mid {q-1}$ and \[q>p+q-1\Big(\sum_{r\, \textnormal{is prime:}\,r\mid{q-1}}\frac{1}{r}\Big),\] there exists a dynamically irreducible set over $\mathbb{F}_q$ with cardinality $p^2$.
    \end{theorem}
To check the dynamical irreducibility of a set of polynomials, we need to verify whether the elements of \eqref{dynamical irr} are $rad(t)$-free. In the following result, we give a bound on the set \eqref{dynamical irr} and it is independent of the size of the set of polynomials.
\begin{theorem}\label{Main4}
    Let $t\geq{3}$ be odd and let \( A=\{f_i=(x-a_i)^t+b_i\in{\mathbb{F}_q[x], 1\leq{i}\leq{k}}\},\) be a dynamically irreducible set of polynomials over a finite field \( \mathbb{F}_q \) of odd characteristic. Assume that $rad(t)\mid q-1$ and \( k \geq 2 \). Suppose $B$ denotes the set $\{f_{j_1}\circ f_{j_2}\circ f_{j_3}\circ\cdots\circ f_{j_r}(b_i):1\leq{j_1,j_2,\ldots,j_r,i}\leq{k} \}$ together with the values of $b_i$. Then, for any prime divisor $r$ of $t$, $\#B=O(q^{1-\alpha_{t,r}})$ where $\alpha_{t,r}=\frac{\log (\frac{r}{r-1})}{2\log {\frac{tr}{r-1}}}$

\end{theorem}

\textbf{Outline of proofs and organization of the paper:}
In Theorem \ref{Main1}, a necessary and sufficient condition is proved for the dynamical irreducibility of a set of polynomials over a finite field. The proof is done by the repeated use of the Capelli lemma and properties of Norm in a finite field. In Theorem \ref{Main2}, an explicit family of dynamically irreducible polynomials in $\mathbb{F}_{p^p}[x]$ is constructed using the criterion proved in Theorem \ref{Main1}. The construction in Theorem \ref{Main2} is motivated by the technique of Heath Brown et. al. \cite{Heath2019}. In Theorem \ref{lower-bound}, the idea of Theorem \ref{Main2} is generalized to construct a dynamically irreducible family of polynomials of certain degrees over any finite field. Moreover, using character sums, the existence of a dynamically irreducible set of polynomials of any degree over any finite field is established. In Theorem \ref{Main4}, using the technique of \cite{Alina} (which is also used in \cite{Jamie}), a bound is given on the cardinality of the set that needs to be tested for the dynamical irreducibility of a family of polynomials. The article is organized as follows. In Section \ref{2}, we recall some known facts that are required to establish our results. In Section \ref{3}, we prove the main results of the article.

\section{preliminary}\label{2}
\begin{defn}
  For \(\alpha \in \mathbb{F}_{q^n}\), the norm of $\alpha$, \(N_{q^n/q}(\alpha)\), over \(\mathbb{F}_q\) is  
\[  
N_{q^n/q}(\alpha) = \alpha \cdot \alpha^q \cdots \alpha^{q^{n-1}} = \alpha^{\frac{q^n-1}{q-1}}.  
\]
\end{defn}
The next result follows from Lemma 3.1 in \cite{Avnish2022} and the fact that \( N_{q^n/q}(\alpha) = 0 \) if and only if \( \alpha = 0 \).
\begin{lemma}
    Let \( t \mid q^n - 1 \), \( \sigma = \gcd(t, q - 1) \), and \( Q_t \) be the largest divisor of \( t \) such that \( \gcd(Q_t, \sigma) = 1 \). Then an element \( \alpha \in \mathbb{F}_{q^n} \) is \( t \)-free if and only if \( \alpha \) is \( Q_t \)-free and \( N_{q^n/q}(\alpha) \) is \( \sigma \)-free in \( \mathbb{F}_q \).
\end{lemma}	
In particular, if $t\mid q-1$, we have the following corollary:
\begin{coro}\label{norm}
    Let \( t \) be an integer dividing \( q - 1 \). Then, \( \alpha \in \mathbb{F}_{q^n} \) is \(t\)-free if and only if \( N_{\mathbb{F}_{q^n}/\mathbb{F}_q}(\alpha) \) is \(t\)-free in \(\mathbb{F}_q\).
\end{coro}
We shall use the following results on the irreducibility of binomials.	

\begin{theorem}[\cite{lidl}, Proposition 3.75]
	Let $e$ be the order of $a\in{\mathbb{F}_q^{*}}$ with respect to multiplication. Then $x^t-a$, for $t\geq{2},$ is irreducible over $\mathbb{F}_q$ if and only if the following conditions are satisfied:
	\begin{itemize}
		\item [1.] If $r$ is a prime divisor of $t$, then $r\mid e$ but $r\nmid\,(q-1)/e.$
		\item [2.] If $4\mid t$, then $q\equiv1\mod{4}$.
	\end{itemize}
\end{theorem}
It is not difficult to observe that the above theorem can be restated as follows:
\begin{prop}\label{bino-irr}
    Let \( t \geq 2 \) be an integer and \( a \in \mathbb{F}_q \). Then, the binomial \( x^t - a \) is irreducible in \( \mathbb{F}_q[x] \) if and only if the following conditions are satisfied:  
(i) \(\text{rad}(t) \mid q - 1\);  
(ii) \( a \) is \(\text{rad}(t)\)-free;  
(iii) \( q \equiv 1 \pmod{4} \) if \( t \equiv 0 \pmod{4} \).
\end{prop}
The next lemma gives a necessary and sufficient condition for the irreducibility of the composition of two polynomials over any field $K$. 
	\begin{lemma}[Capelli's Lemma, \cite{jones2008}]\label{Capelli}
	Suppose $f(x), g(x)$ are polynomials over the field $K$ such that $g(x)$ is irreducible. Then $g(f(x))$ is irreducible over $K$ if and only if $f(x)-\beta$ is irreducible over $K(\beta)$ for every root $\beta\in{\overline{K}}$ of $g(x)$.
\end{lemma}

\begin{lemma}[Lemma 4.1, \cite{Panario2025}]\label{1_{r-free}}

Let $t$ be a prime such that $t\mid q-1$, and let $\chi_t$ be a multiplicative character of $\mathbb{F}_q$ of order $t$. For $\alpha \in \mathbb{F}_q^*$, define
\[
1_{t\text{-free}}(\alpha) := \frac{t-1}{t} \left( 1 - \frac{1}{t-1} \sum_{j=1}^{t-1} \chi_t^j(\alpha) \right).
\]

Then
\[
1_{t\text{-free}}(\alpha) =
\begin{cases}
1, & \text{if } \alpha \text{ is $t$-free},\\[1ex]
0, & \text{if } \alpha \text{ is not $t$-free}.
\end{cases}
\]

\end{lemma}

The following result gives Weil's bound on non-trivial character sums with polynomial arguments:
\begin{theorem}[\cite{lidl}, Theorem 5.41]\label{weil sum}
    Let \( \chi_r \) be a multiplicative character of \( \mathbb{F}_q \) of order \( r > 1 \). Assume \( f \in \mathbb{F}_q[x] \) is a non-constant monic polynomial that is not a \( r \)-th power of a polynomial. If \( e \) is the number of distinct roots of \( f \) in its splitting field over \( \mathbb{F}_q \), then we have $$\Big|\sum_{a\in{\mathbb{F}_q}}\chi_r(f(a))\Big|\leq({e-1})q^{1/2}$$
\end{theorem}
	
	\section{Dynamically irreducibility of polynomials}\label{3}
    We start this section with the following remark and then prove a necessary and sufficient condition for the dynamical irreducibility of a set of certain polynomials over finite fields.
\begin{rem}
    We say that $f_1\circ f_{2}\circ f_{3}\circ\cdots\circ f_{{r}}(x)$ is reducible of minimal length if $f_1\circ f_{2}\circ f_{3}\circ\cdots\circ f_{{r}}(x)$ is reducible but $f_1\circ f_{2}\circ f_{3}\circ\cdots\circ f_{{r-1}}(x)$ is irreducible, and similarly, we say that $f_1\circ f_{2}\circ f_{3}\circ\cdots\circ f_{{r-1}}(b_r)$ is not $t$-free of minimal length if $f_1\circ f_{2}\circ f_{3}\circ\cdots\circ f_{{r-1}}(b_r)$ is not $t$-free but $f_1\circ f_{2}\circ f_{3}\circ\cdots\circ f_{{i}}(b_{i+1})$ is $t$-free for $1\leq{i}\leq{r-2}$.
\end{rem}

% \begin{theorem}\label{Main1}
%     Let $p$ be prime and let $q=p^r$ for some $r\geq{1}$. Assume that $f_i(x)=(x-a_i)^t-b_i$ is an irreducible polynomial over $\mathbb{F}_q$ for each $1\leq{i}\leq{k}$. The set $A=\{f_i(x):1\leq{i}\leq{k}\}$ is dynamically irreducible over $\mathbb{F}_q$ if and only if \begin{equation}\label{dynamical irr}
%         f_{j_1}\circ f_{j_2}\circ f_{j_3}\circ\cdots\circ f_{j_r}(b_i)
%     \end{equation} is rad(t)-free for all $1\leq{j_1,j_2,\ldots,j_r,i}\leq{k}.$
     
% \end{theorem}

\begin{pff}\textbf{of Theorem \ref{Main1}:}

   Let $f_i(x)=(x-a_i)^t+b_i\in{\mathbb{F}_q[x]}$ and $A=\{f_i(x):1\leq{i}\leq{k}\}.$ We prove the negation and show that a polynomial composition in $A$, $f_{j_1}\circ f_{j_2}\circ f_{j_3}\circ\cdots\circ f_{j_r}(x)$, is reducible of minimal length if and only if $f_{j_1}\circ f_{j_2}\circ f_{j_3}\circ\cdots\circ f_{j_{r-1}}(b_{j_r})$ is not $rad(t)-$free of minimal length. Suppose that the composition $f_{j_1}\circ f_{j_2}\circ f_{j_3}\circ\cdots\circ f_{j_r}(x)$ is reducible of minimal length, that is, $f_{j_1}\circ f_{j_2}\circ f_{j_3}\circ\cdots\circ f_{j_r}(x)$ is reducible but $f_{j_1}\circ f_{j_2}\circ f_{j_3}\circ\cdots\circ f_{j_{r-1}}(x)$ is irreducible over $\mathbb{F}_q$. Let $\gamma_1$ be a root of $f_{j_1}(x)$ and let $\gamma_{i+1}$, $1\leq{i}\leq{{r-1}}$, be a root of $f_{j_{i+1}}-\gamma_{i}$. Set $\alpha_1, \alpha_2, \dots,\alpha_r$ as follows:
    \begin{align*}	\alpha_1&=b_{j_1}\in{\mathbb{F}_q},\\
	\alpha_2&=b_{j_2}+\gamma_1\in{\mathbb{F}_q(\gamma_1)},\\ %={\mathbb{F}_{q^{n_1}}}
	\alpha_3&=b_{j_3}+\gamma_2\in{\mathbb{F}_q(\gamma_1,\gamma_{2})},\\%={\mathbb{F}_{q^{n_1n_2}}},\\
		&\cdots&\\
		\alpha_{r}&=b_{j_r}+\gamma_{r-1}\in{\mathbb{F}_q(\gamma_1,\gamma_{2},\ldots,\gamma_{r-1})},\\%={\mathbb{F}_{q^{n_1n_2 \cdots n_{r-1}}}}\,.
		\end{align*}
        In view of Lemma \ref{Capelli} and Proposition \ref{bino-irr}, the irreducibility of $f_{j_1} \circ f_{j_2}\circ\ldots \circ f_{j_{r-1}}$ is equivalent to that $\alpha_{i}$ is $rad(t)$-free in its field of definition for $1\leq{i}\leq{r-1}$. On the other hand, using the same argument, by the reducibility of $f_{j_1}\circ f_{j_2}\circ\cdots \circ f_{j_{r}}$, $\alpha_{r}=b_{j_r}+\gamma_{r-1}\in{\mathbb{F}_q(\gamma_1,\gamma_{2},\cdots,\gamma_{r-1})}$ is not $rad(t)$-free. For $0\leq j<i$, let $$N_{i}^{j}:\mathbb{F}_q(\gamma_1,\gamma_{2},\ldots,\gamma_{i})\to\mathbb{F}_q(\gamma_1,\gamma_{2},\ldots,\gamma_{j})$$ be the norm map. Then, it is not difficult to observe that $$ N_{k}^{0}(\alpha_{k+1})=f_{j_1}\circ  f_{j_2} \cdots \circ f_{j_{k}}(b_{j_{k+1}}),$$ for $1\leq{k}\leq{r-1}.$ By Corollary \ref{norm}, $\alpha_k\in{\mathbb{F}_q(\gamma_1,\gamma_{2},\cdots,\gamma_{k-1})}$ is $rad(t)-$free for $1\leq k\leq{r-1}$ but $\alpha_r\in{\mathbb{F}_q(\gamma_1,\gamma_{2},\cdots,\gamma_{r-1})}$ is not $rad(t)-$free. That is, $f_{j_1}\circ f_{j_2}\circ f_{j_3}\circ\cdots\circ f_{j_{r-1}}(b_{j_r})$ is not $rad(t)-$ free of minimal length, so we are done one way. The converse follows a similar line. Therefore, the set $A$ is dynamically irreducible if and only if $f_{j_1}\circ  f_{j_2} \cdots \circ f_{j_{r}}(b_i)$ is $rad(t)-$free for all $1\leq{j_1,j_2,\ldots,j_r,i\leq k}.$
               
\end{pff}

Next, we construct a dynamically irreducible family of polynomials over the finite field $\mathbb{F}_{p^p}$ where $p$ is a prime.

\begin{pff}\textbf{of Theorem \ref{Main2}:}

    For the dynamical irreducibility of $A$, the irreducibility of $f_{a,b}(x)$ is necessary for each $a,b\in{\mathbb{F}_p}$. The polynomial $f_{a,b}(x)$ is irreducible in $\mathbb{F}_{p^p}[x]$ if and only if (i) $-b-\beta$ is $rad(t)$-free, (ii) $rad(t)\mid p^p-1$ and (iii) if $4\mid t$ then $4\mid p^p-1$. Let $g_{a,b}(x)=(x-a)^t+b\in{\mathbb{F}_p[x]}$. Define $G(x)=g_{a_1,b_1}\circ g_{a_2,b_2}\circ \cdots \circ g_{a_n,b_n}(x)$ and $F(x)=f_{a_1,b_1}\circ f_{a_2,b_2}\circ \cdots \circ f_{a_n,b_n}(x)$ for $n\in{\mathbb{N}}$. It is easy to observe that $F(x)=G(x-\beta)+\beta$. For any $b\in\mathbb{F}_p$, $F(b+\beta)=G(b)+\beta$. Since $G(b)\in{\mathbb{F}_p}$ for any $b\in{\mathbb{F}_p}$, $F(b+\beta)\in{\mathbb{F}_p+\beta}.$ By Theorem \ref{Main1}, the set $A$ is dynamically irreducible if and only if (i) all the elements of $\mathbb{F}_p+\beta$ and $\mathbb{F}_p-\beta$ are $rad(t)$-free, (ii) $rad(t)\mid p^p-1$ and (iii) if $4\mid t$ then $4\mid p^p-1$. The last two conditions are met by the irreducibility of polynomials $f_{a,b}(x)$. By Corollary \ref{norm}, $a\pm\beta\in{\mathbb{F}_p\pm\beta}$ is $rad(t)$-free in $\mathbb{F}_{p^p}$ if and only if $N_{\mathbb{F}_{p^p}/\mathbb{F}_{p}}(a\pm\beta)$ is $rad(t)$-free in $\mathbb{F}_p$. Since, $N_{\mathbb{F}_{p^p}/\mathbb{F}_{p}}(a\pm\beta)=\pm{h}$. Therefore, the set $A$ is dynamically irreducible if and only if (i) $h$ and $-h$ are $rad(t)$-free, (ii) $rad(t)\mid p^p-1$ and (iii) if $4\mid t$ then $4\mid p^p-1$.
\end{pff}

Assume that $\beta$ is same as in Theorem \ref{Main2}. We have the following corollaries to the above result:
\begin{coro}
    Let $t\geq{3}$ be an odd integer, and let $p$ be an odd prime. Assume that $f_{a,b}(x)=(x-a-\beta)^t+b+\beta\in{\mathbb{F}_{p^p}[x]}$. If $rad(t)\mid p-1$, and $h$ is $rad(t)$-free in $\mathbb{F}_p$, then $A=\{f_{a,b}(x):a,b\in{\mathbb{F}_p}\}$ is a dynamically irreducible set in $\mathbb{F}_{p^p}$.
\end{coro}
\begin{pff}
    Since $t$ is odd, there is nothing to prove for condition (iii) of Theorem \ref{Main2}. Since $rad(t)\mid p-1$, it is not difficult to observe that $h$ is $rad(t)$-free if and only if $-h$ is $rad(t)$-free. Thus, the result follows.
\end{pff}

\begin{coro}
    Let $t$ be an odd prime and let $p$ be a prime such that $p\equiv{1} \pmod{t}.$ Let $f_{a,b}(x)=(x-a-\beta)^t+b+\beta\in{\mathbb{F}_{p^p}[x]}$. If $h$ is $rad(t)$-free in $\mathbb{F}_p$, then the set $A=\{f_{a,b}(x):a,b\in{\mathbb{F}_p}\}$ is dynamically irreducible over $\mathbb{F}_{p^p}$.
\end{coro}
\begin{rem}\label{imp_rem}
    As we proved in the previous theorem, for the dynamical irreducibility of the set $\{f_{a,b}(x)=(x-a)^t+b+\beta:a,b\in{\mathbb{F}_p}\}$ over $\mathbb{F}_{p^p}$, it is sufficient to show that every element of $\mathbb{F}_p+\beta$ and $\mathbb{F}_p-\beta$ is $rad(t)$-free in $\mathbb{F}_{p^p}$. We can generalize this idea to an arbitrary finite field $\mathbb{F}_q$, that is, if there exists $\alpha\in\mathbb{F}_q$ such that $\mathbb{F}_p+\alpha$ is $rad(t)$-free in $\mathbb{F}_{q}$ then the set $\{f_{a,b}(x)=(x-a)^t+b+\alpha:a,b\in{\mathbb{F}_p}\}$ is dynamically irreducible over $\mathbb{F}_q$. Note that when $m\geq{3}$ is odd, then each element of $\mathbb{F}_p+\beta$ is $m$-free if and only if each element of $\mathbb{F}_p-\beta$ is $m$-free.
\end{rem}
In Theorem \ref{Main2}, we gave an explicit construction of a dynamically irreducible set over a finite field $\mathbb{F}_{p^p}$. It will be interesting to investigate the existence of a dynamically irreducible set of higher-degree polynomials over an arbitrary finite field. Thus, we prove Theorem \ref{lower-bound}. 
% \begin{theorem}\label{lower-bound}
%     Let $q=p^r$ and let $t=s^k$ be an odd prime power such that $s\mid q-1$. Then, there exists an $\beta\in{\mathbb{F}_q}\setminus \mathbb{F}_p$ such that the set $A=\{(x-a-\beta)^t+b+\beta:a,b\in{\mathbb{F}_p}\}$ is dynamically irreducible. Moreover, for any odd $t\geq{3}$ with $rad(t)\mid {q-1}$ and \[q>p\Big(1+q-1\Big(\sum_{r\, \textnormal{is prime:}\,r\mid{q-1}}\frac{1}{r}\Big)\Big),\] there exists a dynamically irreducible set over $\mathbb{F}_q$ with cardinality $p^2$.
    
%     % polynomials of any odd degree $t$ such that $rad(t)\mid q-1$.
%     % for any $C<(\log{4})^2$ and for any sufficiently large characteristic, there exist infinitely many finite fields such that $$M(q,t)\geq{C(\log{q})^2}.$$
% \end{theorem}

\begin{pff}\textbf{of Theorem \ref{lower-bound}:}

   Our goal is to find a $\beta$ such that every element of $\mathbb{F}_p+\beta$ is $rad(t)$-free. That is, we need to show that \(a+\beta\) is not a $s$-th power in $\mathbb{F}_q$ for all $a\in{\mathbb{F}_p}$. To prove this, we will use the estimates for character sums. Let $\chi_s$ denote a character of order $s$, and consider

\begin{equation}\label{S_p equ.}
    S =
\sum_{\beta \in F_q} \prod_{a \in F_p} (1 - \chi_s(a+\beta)) .
\end{equation}

It is not difficult to observe that if for every $\beta\in{\mathbb{F}_q}$ there exist some $a\in\mathbb{F}_p$ such that $a+\beta$ is a $s$-power then $S=0$. Therefore, we need to show that \(|S| > 0\), so that there must be some \(\beta \in \mathbb{F}_q \setminus \mathbb{F}_p\) such that \(a + \beta\) is not $s$-th power for every \(a\) in $\mathbb{F}_p$.

By expanding the product in \eqref{S_p equ.}, we get

\[
S = q +
\sum_H i(H)\sum_{\beta\in{\mathbb{F}_q}}\chi_s(H(\beta)),
\]

where \(H\) runs over polynomials of the form
\[ H(x)=  \prod_{a \in A} (x + a) \] for all nonempty subsets \(A \subseteq \mathbb{F}_p\), and \(i(H) = (-1)^{\#A}\). We have,
$$|S|\leq q+\left|\sum_H i(H)\sum_{\beta\in{\mathbb{F}_q}}\chi_s(H(\beta))\right|.$$
By Weil’s bound for character sums, we have
\[ \left| \sum_{\beta \in \mathbb{F}_q} \chi_s(H(\beta)) \right| \leq \deg(H)q^{1/2} , \]
which implies that
\[ |S| \geq q - q^{1/2} \sum_{k=1}^{p} \binom{p}{k} = q - q^{1/2}\, p \, 2^{p-1} . \] Therefore, a suitable \(\beta\) exists if \(q \geq p^2{4^p}\). In particular, if \(r \geq 3 + \frac{p \log 4}{\log p}\), then for every characteristic, there is a dynamically irreducible set of polynomials of degree $t$ in $\mathbb{F}_q[x]$ where $q=p^r$. 
\par
Next, assume that $s_1$ and $s_2$ are the only distinct odd primes dividing $t$. Since $rad(t)$ divides $q-1$, we have multiplicative characters $\chi_{s_1},\chi_{s_2}$ of orders $s_1$ and $s_2$, respectively. If there exists a $\alpha\in{\mathbb{F}_q\setminus \mathbb{F}_p}$ such that $a+\alpha \notin {\mathbb{F}_q^{s_i}}$ for each $a\in{\mathbb{F}_p}$ and for $i=1,2$, then each element of $\mathbb{F}_p+\alpha$ is $rad(t)$-free. Therefore, by a similar argument to that in Remark \ref{imp_rem} we get a dynamically irreducible set in $\mathbb{F}_q[x]$. Therefore, our goal is to find a common $\alpha\in{\mathbb{F}_q\setminus \mathbb{F}_p}$ such that no element of $\mathbb{F}_p+\alpha$ belongs to the sets ${\mathbb{F}_q^{s_i}}$ for $i=1,2$. Let $$B_i=\{\alpha\in{\mathbb{F}_q}\mid \,\exists \,a\in{\mathbb{F}_p}\, \textnormal{such that}\,\, \chi_{s_i}(a+\alpha)=1\}.$$ It is sufficient to show that there exists a  $\alpha\in{\mathbb{F}_q\setminus \mathbb{F}_p}$ which does not belong to  $(B_1\cup B_2),$ that is, to show $\mathbb{F}_q\setminus \mathbb{F}_p\ne B_1\cup B_2$. The number of elements in $\mathbb{F}_q$ which are not $s_i$-th power is $1+\frac{q-1}{s_i}$. Therefore,
% For each $a\in{\mathbb{F}_p},$ the number of $\alpha\in{\mathbb{F}_q\setminus \mathbb{F}_p}$ such that $a+\alpha$ is a $s_i$-th power is bounded above by $\frac{q-1}{s_i}$. Since there are $p$ choices for
% $a$,
we have $\mid B_{i}\mid \leq \frac{q-1}{s_i}$. So that 
$$\mid B_1\cup B_2\mid\leq (q-1)\Big(\frac{1}{s_1}+\frac{1}{s_2}\Big).$$
For the existence of one $\alpha$ in $\mathbb{F}_q\setminus \mathbb{F}_p$, it is sufficient that $(q-1)\Big(\frac{1}{s_1}+\frac{1}{s_2}\Big)<q-p$. That is, if $q>p+q-1(\frac{1}{s_1}+\frac{1}{s_2})$ then there exist a $\alpha$ such that $a+\alpha$ is  $rad(t)$-free. More generally if $t$ be any odd integer and $p_1,p_2,\cdots, p_n$ are the distinct prime divisors of $t$ then there exist a $\alpha\in{\mathbb{F}_q\setminus {\mathbb{F}_p}}$ such that $a+\alpha$ is $rad(t)$-free for each $a\in{\mathbb{F}_p}$ if $q>p+q-1(\sum_{i=1}^{n}\frac{1}{p_i})$ and thus we have a dynamically irreducible set of polynomials of degree $t$.
\end{pff}
\begin{rem}
    Let $M(q,t)$ be the cardinality of the maximal set of dynamically irreducible polynomials of degree $t$ in $\mathbb{F}_q[x]$. For each characteristic $p$, $M(q,s^n)\geq{p^2}$ for any odd prime $s\mid q-1$ for all $q$ except finitely many. Moreover, for certain choices of $q$, we have $M(q,t)\geq{p^2}$ for any odd $t\geq{3}$ such that $rad(t)\mid q-1$.
\end{rem}
The idea for the proof of the following result is motivated from the work on decomposable polynomials \cite{Von zur} and a similar proof may be found in the literature, but the proof is written here for completeness.
\begin{lemma}\label{composition of two}
    Let $f_1,f_2$ be two distinct monic irreducible polynomials of degree $t\geq{2}$ of the form $(x-a)^t+b$ in $\mathbb{F}_q[x]$ where $q$ is odd. If 
    \begin{equation}\label{equ composition}
         f_{i_1} \circ \cdots \circ f_{i_n} = f_{j_1} \circ \cdots \circ f_{j_m} 
    \end{equation}
with \( i_1, \ldots, i_n, j_1, \ldots, j_m \in \{1, 2\} \). Then \( m = n \) and \( i_k = j_k \) for every index \( k\).
\end{lemma}
\begin{pff}
    Since Equation \eqref{equ composition} holds, the degrees of polynomials on both sides are the same, and thus $m=n$. Assume to the contrary, $n$ is minimal such that the order of composition of the polynomials $f_1$ and $f_2$ does not preserve, that is, a nontrivial relation as in Equation \eqref{equ composition} holds. Then, $f_i\circ F=f_j\circ G$ where $f_i\ne f_j$ or $f_i= f_j$ but $F\ne G$. Let $f_i(x)=(x-a)^t+b$ and $f_j(x)=(x-c)^t+d$. Then,
    \begin{align*}
        b-d&=(G(x)-c)^t-(F(x)-a)^t\\
        &=(G(x)-F(x)+a-c)\Big((G(x)-c)^{t-1}+\cdots+(F(x)-a)^{t-1}\Big)
    \end{align*}
    Since $F(x)$ and $G(x)$ are monic and $\mathbb{F}_q$ has odd characteristic, $(G(x)-c)^{t-1}+\cdots+(F(x)-a)^{t-1}$ has positive degree. Therefore, we have $b=d$ and $G(x)-F(x)+a-c=0$. If $a=c$, we have $f_i=f_j$ and $F=G$, which is a contradiction to the assumption. Therefore, $a\ne c$ so that we have $f_1(x)=(x-a)^t+b$, $f_2(x)=(x-c)^t+b$ in some order, and $F\ne G$. Without loss of generality, we may assume that $F=f_r\circ F_1$ and $G=f_s\circ G_1$ where $f_r(x)=(x-a)^t+b$, $f_s(x)=(x-c)^t+b$. Then,
    \begin{align*}
0 \ne a - c =& F(X) - G(X) = (F_1(X) - a)^t - (G_1(X) - c)^t \\
=& \Big((F_1(X)-a)^{t-1} +\cdots (G_1(X)-c)^{t-1}\Big)\Big(F_1(X) - G_1(X) - a+ c\Big).
    \end{align*}
    Since the polynomial $(F_1(X)-a)^{t-1} +\cdots +(G_1(X)-c)^{t-1}$ has a positive degree, we have a contradiction. 
\end{pff}

For the dynamical irreducibility of the set $\{f_i(x)=(x-a_i)^t+b_i:1\leq{i}\leq{k}\}$, we need to verify whether the elements in \eqref{dynamical irr} are $rad(t)$-free in $\mathbb{F}_q$. In the following result, we give a bound on the cardinality of the set \eqref{dynamical irr}, and some construction is inspired by the ideas of \cite{Panario2025}. Moreover, this bound is independent of the size of the set of polynomials. 
% \begin{theorem}
%     Let $t\geq{3}$ be odd and let \( A=\{f_i=(x-a_i)^t-b_i\in{\mathbb{F}_q[x], 1\leq{i}\leq{k}}\},\) be a dynamically irreducible set of polynomials over a finite field \( \mathbb{F}_q \) of odd characteristic. Assume that $rad(t)\mid q-1$ and \( k \geq 2 \). Suppose $B$ denotes the set $\{f_{j_1}\circ f_{j_2}\circ f_{j_3}\circ\cdots\circ f_{j_r}(b_i):1\leq{j_1,j_2,\ldots,j_r,i}\leq{k} \}$ together with the values of $b_i$. Then, for any prime divisor $r$ of $t$, $\#B=O(q^{1-\alpha_{t,r}})$ where $\alpha_{t,r}=\frac{\log (\frac{r}{r-1})}{2\log {\frac{tr}{r-1}}}$

% \end{theorem}
\begin{pff}\textbf{of Theorem \ref{Main4}:}

    Note that $B\subset{\mathbb{F}_q}$. Due to the dynamical irreducibility of the set $A$, Theorem \ref{Main1} implies that $f_{j_1}\circ f_{j_2}\circ f_{j_3}\circ\cdots\circ f_{j_r}(\beta)$ is $rad(t)$-free in $\mathbb{F}_q$ for any $\beta\in{B}$. Let $f_1,f_2$ be fixed polynomials in $A$, and let $N$ be a fixed natural number. Assume that $\mathfrak{F}$ denotes the set of distinct polynomials given as $N$ time composition of $f_1$ and $f_2$. By Lemma \ref{composition of two}, $|\mathfrak{F}|=2^N$. In particular, $\mathfrak{F}$ is the set of irreducible polynomials of degree $t^N$ given as compositions of $f_1$ and $f_2$. For a fixed prime $r$ dividing $t$ we consider the following set:
    $$\tau_{q,r}(N)=\{a\in{\mathbb{F}_q}:{1_{rf}}(F(a))=1, \forall \, F\in{\mathfrak{F}}\}.$$ By Theorem \ref{Main1}, $B\subseteq{\tau_{q,r}(N)}$ and so $|B|\leq{|\tau_{q,r}(N)|}$. The following equality is easy to observe: $$|\tau_{q,r}(N)|=\sum_{a\in{\mathbb{F}_q}}\Big(\prod_{s=1}^{2^N}1_{r-free}\big(F_s(a)\big)\Big)$$
    Since $F_s(a)\in{\mathbb{F}}^{*}_{q},$ by Lemma \ref{1_{r-free}} we have 
    \begin{equation}\label{eq1}
      |\tau_{q,r}(N)|= \sum_{a\in{\mathbb{F}_q}}\Big(\prod_{s=1}^{2^N}\frac{r-1}{r}\Big(1-\frac{1}{r-1}\Big(\sum_{j=1}^{r-1}\chi_r^j(F_s(a))\Big)\Big)\Big).  
    \end{equation}
    Now, set $\omega(N)=\{0,1,2,\ldots,r-1\}^{2^N}$ and for each $0\leq{u}\leq{2^N}$, set $\omega(N)=\{c\in \omega(N): c\, \,\textnormal{has exactly}\, u \,\,\textnormal{vanishing coordinates}\}.$
    For $c=(c_1,c_2,\ldots,c_{2^N})\in \omega(N)$, let 
    $$G_{F,c}=\sum_{a\in{\mathbb{F}_q}}\chi_r\Big(\prod_{s=1}^{2^N}F_s(a)^{c_s}\Big).$$
    By expanding the product in \eqref{eq1}, we get
    \begin{equation}\label{eq2}
        |\tau_{q,r}(N)|=\Big(\frac{r-1}{r}\Big)^{2^N}\Bigg(q+\sum_{u=1}^{2^N}\Big(\frac{-1}{r-1}\Big)^u\sum_{c\in \omega(N)}G_{F,c}\Bigg).
    \end{equation}
    Since the polynomials in $\mathfrak{F}$ are all irreducible in $\mathbb{F}_q[x]$, no two polynomials share a common root. Thus, the polynomial $\prod_{s=1}^{2^N}F_s(a)^{c_s}$ is not a $r$-th power in the algebraic closure of $\mathbb{F}_q$. Additionally, the polynomial $\prod_{s=1}^{2^N}F_s(a)^{c_s}$ has at most $(2t)^N$ distinct roots for each $1\leq{u}\leq{2^N}$ and each $c=(c_1,c_2,\ldots,c_{2^N})\in{\omega_u(N)}$. By Weil's bound for character sums, we have 
    \begin{equation}\label{eq3}
        |G_{F,c}|\leq{(2t)^N}q^{1/2}.
    \end{equation}
    From Equations \eqref{eq2},\eqref{eq3}, we have 
    \begin{align*}
        |\tau_{q,r}(N)|&\leq \Big(\frac{r-1}{r}\Big)^{2^N}q+{\Big(\frac{r-1}{r}\Big)^{2^N}\sum_{u=0}^{2^N}\binom{2^N}{u}\Big(\frac{1}{r-1}\Big)^{u}}(2t)^Nq^{1/2}\\
        &=\Big(\frac{r-1}{r}\Big)^{2^N}q+\Big(\frac{r-1}{r}\Big)^{2^N}\Big(1-\frac{r}{r-1}\Big)^{2^N}(2t)^Nq^{1/2}\\
        &=\Big(\frac{r-1}{r}\Big)^{2^N}q+(2t)^Nq^{1/2},
    \end{align*}
%     then
%     \[
% \Big(\frac{r-1}{r}\Big)^{2^N} q-(2t)^N q^{1/2} \leq \#\tau_{q,r}(N) \leq \Big(\frac{r-1}{r}\Big)^{2^N}q+(2t)^N  q^{1/2}.
% \]
Choose $N$ such that $$2^N\leq{\frac{1}{2} \log_{\frac{tr}{r-1}}q}<2^N+1$$ so that $\Big(\frac{r-1}{r}\Big)^{2^N}q\leq (t^{2^N}q^{1/2}).$ Therefore, we have $\# \tau_{q,r}(N) = O(t^{2^N}{q^{1/2}}).$
$$ t^{2^N}q^{1/2} \leq t^{\frac{1}{2} \Big(\log_{\frac{tr}{r-1}}q\Big)}q^{1/2}=q^{\frac{1}{2}\Big(1+\frac{\log t}{ \log{\frac{tr}{r-1}}}\Big)}=q^{1-\alpha_{t,r}} 
$$ where $\alpha_{t,r}=\frac{\log (\frac{r}{r-1})}{2\log {\frac{tr}{r-1}}}$.
\end{pff}
% Note that the explicit construction in Theorem \ref{Main2} implies that $M(p^p)\geq{p^2}$. Some of the constructions of the proof are inspired by the arguments of Heath brown \cite{HB}.
% \par
% Due to Theorem \ref{lower-bound}, we have a dynamically irreducible set of polynomials of prime power degree over infinitely many finite fields with a sufficiently large characteristic. Next, we want to ask when we can have a dynamically irreducible set of polynomials of degree $t$ where $t$ is not a prime power.  

% In Theorem \ref{Main1}, we have a necessary and sufficient condition for the dynamical irreducibility of a set of binomials over a finite field. Although we do not have a necessary and sufficient condition over an arbitrary field, we have the following sufficient condition for the dynamical irreducibility of binomials over an arbitrary field. 
% \begin{theorem}
%     Let $t\geq{1}$, $r\in{\mathbb{N}}$, and let $K$ be a field. Assume that the polynomials $f_i(x)=(x-a_i)^t-b_i\in{K[x]}$ are irreducible over $K$ for $1\leq{i}\leq{r}$. The set $A=\{f_i(x):1\leq{i}\leq{r}\}$ is dynamically irreducible over $K$ if $f_{j_1}\circ f_{j_2}\circ f_{j_3}\circ\cdots\circ f_{j_s}(b_i)$ is not a $p$-th power for each $1\leq{j_1,j_2,\ldots,j_s,i}\leq{r}$ for any prime divisor $p$ of $t$. 
% \end{theorem}
% \begin{pff}
%     The proof is on the similar lines as of the proof of Theorem \ref{Main1}. By Irr them, Capellis lemma, norm remark the proof follows.
% \end{pff}
% \begin{exmp}
    
% \end{exmp}
% We have the following example of a dynamical irreducible set over $\mathbb{Q}$:

\textbf{Future direction:}\\
In the literature, we have explicit examples of a dynamically irreducible set of polynomials over a very restricted family of finite fields. One can try to construct an explicit family of dynamically irreducible polynomials over any finite field.  \\
The existing literature provides only discussion and explicit examples of dynamically irreducible families of polynomials in which all polynomials have the same degree. A natural subsequent question is whether one can construct an explicit family of dynamically irreducible polynomials without imposing the restriction that the degrees be equal.
\\
\\

% \textbf{Declarations:}
% \begin{itemize}
% 	\item \textbf{Ethical Approval and Consent to participate:} Not applicable.
% 	\item \textbf{Consent for publication:} Not applicable.
% 	\item \textbf{Availability for supporting data:} Not applicable.
% 	\item \textbf{Competing interests:} The authors declare that they have no competing interests.
% 	\item \textbf{Funding:} Not applicable.
% 	\item \textbf{Author's contributions:} All authors contributed equally to this work. All authors read and approved the final manuscript.
% 	\item \textbf{Acknowledgments:} Not applicable.
% \end{itemize}

	\bibliographystyle{plain}

\end{document}